\documentclass[12 pt]{amsart}
\usepackage{amscd,amsfonts,amssymb,amsmath}
\usepackage{hyperref}
\usepackage[cp1251]{inputenc}
\usepackage[T2A]{fontenc}

\usepackage{epsfig}
\newtheorem{theorem}{Theorem}[section]
\newtheorem{corollary}[theorem]{Corollary}
\newtheorem{lemma}[theorem]{Lemma}
\newtheorem{proposition}[theorem]{Proposition}
\theoremstyle{definition}
\newtheorem{conjecture}[theorem]{Conjecture}
\newtheorem{question}[theorem]{Question}

\newtheorem{example}[theorem]{Example}

\numberwithin{equation}{subsection}
\newtheorem*{ack}{Acknowledgement}
\usepackage[all,cmtip]{xy}
\usepackage{xcolor}
\usepackage{xfrac}
\usepackage{pb-diagram}

\newcommand{\Aut}{\operatorname{Aut}}

\def\bea{\begin{eqnarray}}
\def\eea{\end{eqnarray}}

\DeclareMathOperator{\Set}{Set}
\DeclareMathOperator{\Sym}{Sym}

\begin{document}
\title[Cayley graphs]{
Cayley graphs and their growth functions for  multivalued groups
}
\author{Valeriy G. Bardakov}
\author{Tatyana A. Kozlovskaya}
\author{Matvei N. Zonov}

\date{\today}
\address{Sobolev Institute of Mathematics, 4 Acad. Koptyug avenue, Novosibirsk, 630090, Russia.}
\address{Novosibirsk State Agrarian University, Dobrolyubova street, 160, Novosibirsk, 630039, Russia.}
\address{Regional Scientific and Educational Mathematical Center of Tomsk State University,
36 Lenin Ave., 14, 634050, Tomsk, Russia.}
\email{bardakov@math.nsc.ru}

\address{Regional Scientific and Educational Mathematical Center of Tomsk State University, 36 Lenin Ave., 14, 634050, Tomsk, Russia.}
\email{t.kozlovskaya@math.tsu.ru}

\address{Sobolev Institute of Mathematics, 4 Acad. Koptyug avenue, Novosibirsk, 630090, Russia.}
\address{Regional Scientific and Educational Mathematical Center of Tomsk State University,
36 Lenin Ave., 14, 634050, Tomsk, Russia.}
\email{mat.zonov@yandex.ru}

\subjclass[2000]{Primary 20N20; Secondary 16S34, 05E30}
\keywords{Multi-set, multivalued group, Cayley graphs,  $n$-valued dynamic, growth function}

\begin{abstract} We define the Cayley graph  and its  growth function for  multivalued groups. We prove that if we change a finite set of generators  of multivalued  group, or change the starting point, we get an equivalent growth function. We prove that if we take a virtually nilpotent group and construct a coset group with respect a finite group of authomorphisms, then this multivalued group has a polynomial growth.  
Also, we find a connection between this growth function  and growth function of multivalued dynamics. It particular, it is obtained upper and lower bounds  on growth functions of multivalued dynamics. We give a particular  answer to a question of Buchstaber on polynomial growth of dynamics and a question of Buchstaber and Vesnin on growth functions of cyclically presented multivalued groups.
\end{abstract}
\maketitle

\section{Introduction}

 In this paper by 
a groupoid we mean an algebraic system with one binary operation. A groupoid system is an algebraic system with a family of binary operation. A groupoid system 
is an example of algebraic systems  with a set of one type algebraic operations. In \cite{K} algebraic systems of this type was called homogeneous. 
For example,  a brace or a skew brace  is a set with two group operations, which satisfy appropriate axiom (\cite{Rump}, \cite{GV}). A generalization of  brace is a brace systems \cite{BNY} that is   a set with a family of group operations with some consistency conditions.

An example of a semigroup system 
is a dimonoid which was  introduced by J.-L.~Loday~\cite{Loday2}. A dimonoid is a set with two semigroup operations which are connected by a set of axioms. 
Structural properties of free dimonoids was studied by A.~V.~Zhuchok in \cite{Zhu}. 
T.~Pirashvili \cite{Pir} considered sets with two associative binary operations (without added connections between these operations) and introduced the notion of a duplex. 
 Dimonoids in the sense of J.-L.~Loday are examples of duplexes.

An algebraic system  $\mathcal{G} = (G, *_i, i \in I)$, where $(G, *_i)$ is a (semi-)group  for any $i \in I$ is said to be a 
(semi-)group system. 
A semigroup system with two operations is called a duplex. We call $\mathcal{G}$ by $I$-multi-(semi-)group if the operations are related by the axiom of mixed associativity,
$$
(a *_i b) *_j c = a *_i (b *_j c),~~ a, b, c \in G,~~i, j \in I. 
$$
A multi-semigroup with $n$ operations is called by $n$-tuple semigroup in \cite{Kor}.
In \cite{K} it were defined several homogeneous algebraic systems which are generalization of  multi-(semi)-group.

Another type of algebraic systems suggested in 1971, V.~M.~Buchstaber and S.~P.~Novi\-kov~\cite{BN}. They  gave a notion of $n$-valued group in which the product of each pair of elements is an $n$-multi-set, the set of $n$ elements with multiplicities (see Section \ref{prel} for the definition).   In this case we can not define $n$-valued group as a group system \cite{BKT}.
An appropriate survey on $n$--valued groups and its applications can be found in \cite{B}. Also, in this paper was suggested three open questions.
In  \cite{BK} it  was given particular answer to the first two questions and  a complete answer to the third question. 

In \cite{BKT}  we  investigate connections between $n$-multi-groups  and $n$-valued groups.
If we have a group system $\mathcal{G} = (G, *_i, i \in I)$, where $|I| = n$, we can  define a $n$-valued multiplication 
$$
a * b = [a *_1 b, a *_2 b, \ldots,  a *_n b],~~a, b \in G,
$$ 
and study conditions under which the  algebraic system $(G, *)$ is a $n$-valued group.
We prove  that if in the  groups 
$(G, *_i)$, $i = 1, 2, \ldots, n$ all units are equal  and $(G, *)$ is an $n$-valued group, then $*_i = *_j$ for all $1 \leq i, j \leq n$. It means, in particular, that in some sense multi-groups and multivalued groups are principally different algebraic systems.

Also, in \cite{BKT}
 we introduce and study $n$-valued racks and quandles.  We suggest some constructions for them. In particular, a coset quandle, which is similar to
the coset group in the category of groups. The second construction is new,  we prove that if a set $X$ is equipped with $n$ quandle operations $*_i$ such that
\begin{equation} \label{ass}
(x *_i y) *_j z= (x *_j z) *_i (y *_j z),~~~x, y, z \in X,
\end{equation}
then an $n$-valued multiplication 
$$
x * y = [x *_1 y, x *_2 y, \ldots, x *_n y]
$$
defines  an $n$-valued quandle structure on $X$. We  call a set $X$ with a set of quandle operations $(X, *_1, *_2, \ldots, *_n)$ satisfying axioms of mixed associativity (\ref{ass}) a $n$-multi-quandle. Such axioms arise in \cite{BF} in studying multiplications of quandle structures. 
In this paper were considered quandle systems $\mathcal{Q} = (Q, *_i, i \in I)$, where  $(Q, *_i)$ is a quandle for any $i \in I$, and defined a multiplication $*_i *_j$ by the rule
$$
p (*_i *_j) q = (p *_i q) *_j q,~~p, q \in Q.
$$
Generically  $(Q, *_i *_j)$ is not a quandle, but if one imposes additionally
\begin{equation} \label{quand}
(x *_i y) *_j z= (x *_j z) *_i (y *_j z),~~(x *_j y) *_i z= (x *_i z) *_j (y *_i z),~~~x, y, z \in Q,
\end{equation}
then  $(Q, *_i *_j)$ and  $(Q, *_j *_i)$ are quandles. The term $n$-multi-quandle was introduced by V.~G.~Turaev \cite{T}  and gave them a  topological interpretation.

In the present  paper we introduce a Cayley graph for $n$-valued group, give a definition of growth function, which  generalizes the definition of the growth function for $n$-valued dynamic. Using  Gromov's theorem \cite{G} on  polynomial growth of virtually nilpotent groups  (recall that a group is called by virtually nilpotent if it contains a nilpotent subgroup of finite index), we prove that some  $n$-valued groups which constructed by virtually nilpotent groups
have  a polynomial growth. 
Also, we study  growth functions for $n$-valued dynamics that is a particular case of growth function for $n$-valued group. A growth function for 
$n$-valued group great or equal a growth function of any $n$-valued dynamic in this group. In particular, in $n$-valued group of polynomial growth  
any $n$-valued dynamic has a polynomial growth. 
We study dynamics in some 2-valued groups and prove that they have polynomial growth.

The paper is organized  as follows. In Section \ref{prel} we recall a definition of a multivalued group, recall two constructions of multivalued groups: coset groups and double coset groups. 
In section \ref{Cayl} we give definition of the Cayley graph and its growth function for a multivalued group. We prove that the growth function of a Cayley graph depends on the centre  in which we consider the balls. On the other side, we prove (see Proposition  \ref{asimpt}) that if we go from one finite set of generators  to other finite set of generators and from
one center of balls to other center, we get an equivalent growth function.  Theorem~\ref{Gr} says, that if  $G$ is a finitely generated virtually nilpotent group and  $A$ is  a finite subgroup of $\Aut(G)$ of cardinality $n$, then the  $n$-valued coset group $X = (G, A)$   has a polynomial growth.

In section \ref{dyn} we show that for $n$-valued group $X$ the growth function of any $n$--valued dynamic \cite{Bu} is less or equal to the growth function of the Cayley graph of $X$. Then from Theorem \ref{Gr} follows  \label{poldyn} Proposition \ref{poldyn} which says that 
if $X = (G, A)$ is a $n$-valued coset group, then any $n$-valued  dynamic $T_z \colon X \to (X)^n$ which is defined by an element   $z \in X$ has the  polynomial growth function $\xi_y(r)$ in  any $y \in X$.
This proposition  gives   a particular  answer to a question of Buchstaber \cite{Bu}. 
By analogy to this result one can suggest that if $G$ is a group of exponential growth, then any $n$-valued dynamic 
 on the coset group $X = (G, A)$  has an exponential growth. The Example \ref{ex} shows that it is does not hold.

Theorem \ref{Th growth function bounds}
says that in the coset group $X=(G,A)$ for the element $z=gA$ and for any $y\in X$ the growth function $\xi_y(r),$ which is defined by the dynamic $T_z,$ satisfies the inequality:
$$
\frac{1}{n}|S^+(e, r)|\leq\xi_y(r)\leq |B^+(e, r)|,
$$ 
where $S^+(e, r)$ and $B^+(e, r)$ are the sphere and the ball respectively  in the monoid $M=\bigl\langle\{ g a  \mid a\in A\}\bigr\rangle\subset G.$
Using this result, we give an answer on a question of V.~Buchstaber and A.~Vesnin~\cite{BV} on growth functions of cyclically presented multivalued groups.
In Corollary  \ref{cor4.4} we prove that there are cyclically presented groups with dynamics of   exponential growth and cyclically presented groups with cyclic dynamics of a polynomial growth.

The last Theorem \ref{t4.8} gives a  result on growth functions of dynamics in $2$-valued groups. It says that  if 
$X$ is a $2$-valued group such that $inv(a)=a$ for any $a\in X$, then for any $x\in X$ the growth function $\xi_x(r),$ which is defined by the dynamic $T_x,$ satisfies the inequality $\xi_x(r)\leq r(r+1).$ In particular, has a polynomial growth.

\bigskip


\section{Prelimanary results} \label{prel}

\subsection{Multivalued groups}
Recall definitions and some facts from the theory of multivalued groups (see, for example, \cite{B}).

Let $X$ be a non-empty set. An $n$-{\it valued multiplication} on $X$ is a map
$$
\mu \colon X \times X \to (X)^n = \Sym^nX,~~\mu(x, y) = x * y = [z_1, z_2, \ldots, z_n],~~z_k = (x*y)_k,
$$
where $(X)^n = \Sym^nX$ is the $n$-th symmetric power of $X$, that is the quotient $X^n / S_n$ of the Cartesian power $X^n$ under the action of $S_n$ by permutations of components. 
The next axioms are  natural generalizations of the classical axioms of group multiplication.

{\it Associativity}. The $n^2$-multi-sets:
$$
[x * (y * z)_1, x * (y * z)_2,  \ldots, x * (y * z)_n], ~~[ (x * y)_1 * z, (x * y)_2 * z,  \ldots, (x * y)_n * z]
$$
are equal for all $x, y, z \in X$.

{\it Unit}. An element $e \in X$ such that 
$$
e * x = x * e = [x, x, \ldots, x]
$$
 for all $x \in X$.

{\it Inverse}. A map $inv \colon X \to X$ such that 
$$
e \in inv(x) * x~\mbox{and}~e \in x * inv(x)
$$
 for all $x \in X$. Further, for simplicity, we will write $\bar{x}$ instead $ inv(x)$.

The map $\mu$ defines $n$-{\it valued group structure} $\mathcal{X} = (X, \mu, e, inv)$ on $X$ if it is associative, has a unit and an inverse. An $n$-valued group structure on $X$ is commutative if $x * y = y * x$.

\begin{example}
Let us recall two constructions of  $n$-valued groups (see \cite{B}).
 Let $(G; \cdot)$ be a group, $A \leq \Aut(G)$ and $|A| = n$. The set of orbits  $X = G / A$  can be equipped with an $n$-multiplication $\mu \colon X \times X \to (X)^n$ defined by the rule
$$
\mu(x, y) = \pi (\pi^{-1}(x) \cdot \pi^{-1}(y)),
$$
where $\pi \colon G \to X$ is the canonical projection. Then $X$ with multiplication $\mu$ is  a $n$-valued group with the unit  $e_X = \pi(e_G)$ and the inverse $inv(x)$ of $x \in X$ is $\pi (( \pi^{-1}(x))^{-1})$. This $n$-valued group is called the {\it coset group} of the pair $(G, A)$.

The next construction of an $n$-valued group is  called a {\it double coset group}  of a pair $(G, H)$,
 where $G$ is a group and $H$ is its subgroup of cardinality  $n$. Denote by $X$  the space of double coset classes $H \backslash G / H$.
One could define the  $n$-valued multiplication $\mu \colon X \times X \to (X)^n$
by the formula
$$
\mu(x, y) = \{H g_1 H\} * \{H g_2 H\} = [\{H g_1 h g_2 H,~~h \in H],
$$
 and the inverse element $inv_X(x) = \{H g^{-1} H\}$, where 
$x =\{H g H\}$.
\end{example}


\subsection{Cayley graphs of groups and growth  functions}
In this section, we recall some well-known facts from geometric group theory (see, for example, \cite{LS}).

An action of a group $G$ on a set $X$ is a function
$$
G \times X \to X,
$$
where the image of $(g, x)$ is written $g \cdot x$ and 

1) $1 \cdot x = x$ for all $x \in X$;

2) $g \cdot (h \cdot x) = (g h) \cdot x$ for all $g, h \in G$ and $x \in X$.


If  $G$ acts  on a set $X$, then the orbit of $x \in X$ under $G$ is
the set
$$
G \cdot x = \{ g \cdot x~|~ g \in  G \}.
$$
 If no nontrivial element of $G$ fixes any $x \in X$, one says  that the action of $G$ on
$X$ is free.
Any group $G$ acts on itself  by left multiplication: $g\cdot h = g h$. This action is transitive and free.
 
 

Let $G = \langle S \rangle$ be a group with a generating  set $S = S^{-1}$.  The Cayley graph for $G$ with respect to $S$ is a directed, labeled graph $\Gamma_{G,S}$ given as follows: the vertex set is $G$, there is a directed edge from $g$ to $g s$ for every $g \in  G$ and $s \in S$, and we
label the edge from $g$ to $gs$ by the element $s$ of $S$.
Each element of a group $G$ gives rise to an automorphism of the Cayley graph $\Gamma_{G,S}$, as follows. The automorphism $\Phi_g$
associated to $g \in  G$ is given on the vertices by:
$$
\Phi_g(v) = gv.
$$
This is just the action of $G$ on itself by left multiplication.

Using the Cayley graphs we can consider groups as metric spaces.
A closed ball of radius $r \geq 0$ centered at an element $x \in G$ is the set
$$
B(x, r) = \{ y \in G~|~d(x, y) \leq  r \},
$$
where $d(x, y)$ is the length of a shorter way from $x$ to $y$, or $d(x, y) = d(e, x^{-1} y) = l(x^{-1} y)$
is the distance from the identity to $ x^{-1} y$ (in other words, the length of the shortest word representing the element $ x^{-1} y$). Thus, the ball $B(e, r)$ of radius $r$ centered at the identity~$e$ is the set of group
elements whose word length is less than or equal to $r$. There is a bijection $B(e, r) \to B(x, r)$ for any $x\in G$, which is defined by the left multiplication on~$x$. 
The sphere $S(e,r)$ of radius $r$
is the set of all group elements whose word length is exactly $r$.

Two nondecreasing positive functions $f (r)$ and $g(r)$ have the same growth (in this case we write $f(r) \sim g(r)$) if there is a constant $C > 0$ so that
$$
f (r/C) \leq g(r) \leq f (Cr)
$$
for all $r \geq 0$. Then we say that a function $f (r)$ has polynomial growth if there are two positive numbers $k$ and $C$ such that for all $r \geq 1$, one has 
$f(r) \leq C r^k$. One says that a group $G$ has polynomial growth if the function $|B(e, r)|$ has polynomial growth.
The famous Gromov’s  theorem \cite{G} says that a finitely generated
group has polynomial growth if and only if it has a subgroup of finite index that
is nilpotent.

\bigskip


\section{Cayley graph and growth functions for multivalued groups} \label{Cayl}

\subsection{Cayley graph }
 Suppose that $X$ is an $n$-valued group with a finite generating set~$S$ that means that any element of $X$ lies in some product elements of $S$.  The Cayley graph $\Gamma_{X,S} = \Gamma(V, E)$
is an oriented graph with the set of vertices $V = X$ and two vertices $u, v \in V$ are connected by an oriented edge $f \in E$ with a label $s_k$ if for some $s \in S$ and $k \in \{1, 2, \ldots, n \}$ we have $v = (u* s)_k$, where 
$$
u * s = [(u * s)_1, (u * s)_2, \ldots, (u * s)_n].
$$
In this case we say that $u$ is the   initial  vertex of $f$ and $v$ is the terminal vertex of $f$. In fact, for any $s \in S$ and for any $u \in X$ there exists $n$ edges which starting in $u$ and correspond to multiplication on $s$. Of course, some vertices in the multi-set $ [(u * s)_1, (u * s)_2, \ldots, (u * s)_n]$ can be equal. Hence in $\Gamma_{X,S}$ there exists  different edges  which connect one pair of vertices.

For a set $X$ consider a map 
$$
\Set\colon \bigcup\limits_{i=0}^\infty \Sym^i(X)\rightarrow 2^X,
$$
 that maps a $i$-multi-set $M\in \Sym^i(X)$ into the subset of $X$ consisting of all different elements of~$M.$ The number $|\Set(M)|$ can be thought of as a number of different points in $M$.
Note that if $X$ is an $n$-valued group and $A$ is a nonempty subset of $X$ of finite cardinality $k,$ then $A$ can be viewed as an element of $\Sym^k(X),$ and the $n$-valued group operation can be applied to $A$.

\begin{question}
 Let $X$ be a finitely generated $n$-valued group. Is it true that the Cayley graph of $X$ is locally finite?


\end{question}


\subsection{Growth function} We can assume that any edge in the Cayley graph has the length~1 and we can define a distant between vertices $u, v \in V$ as the length of shorter oriented way which connect $u$ with $v$. Remark that by this definition the distant from $u$ to $v$ can be different from the distant from $v$ to $u$. Hence, we are not considering the Cayley graph $\Gamma = \Gamma_{X,S}$ as a metric space, but we can still define 
the {\it length of an element} $x \in X$ as a  minimal $m$ such that $x \in b_1 * b_2 * \ldots * b_m$ for some $b_i \in S$. We will  write $m = l_S(x)$ or simply  $m = l(x)$. Also, we can define a ball of radius $r \geq 0$ with the center in $x\in X$
as the set
$$
B(x, r) = \bigl\{ y\in X~|~\exists m\in \mathbb{N}, m\leq r, \exists s_{i_1},\ldots s_{i_m}\in S, y\in\Set\bigl(x*s_{i_1}*s_{i_2}*,\ldots,*s_{i_m}\bigr) \bigr\},
$$
and a sphere 
$$
S(x, r) = B(x,r)\setminus B(x,r-1).
$$
Then the {\it growth function} of $X$ with respect to  $x \in X$ is a function $\gamma_{\Gamma,x} \colon \mathbb{N} \to \mathbb{N}$, $\gamma_{\Gamma,x}(r)  = |B(x, r)|$. In this definition $r$ can be arbitrary non-negative real  number, but it is evident that $B(x, r) = B(x, [r])$, where $[r]$ is the integer part of $r$. Remark that the growth function depends not only on  $n$-valued group $X$ but also on the generating set $S$. Therefore, it is more accurate to refer to the growth function of the Cayley graph. but later we show (see Proposition \ref{asimpt}) that growth functions with respect to finite generating sets are equivalent and we can say of the growth function of the $n$-valued group $X$.

In the case of ordinary group the left multiplication on some element is an automorphism of Cayley graph and the cardinality $|B(x, r)|$ does not depend on $x$. 
This is not the case for multivalued groups, as the following example shows.

\begin{figure}[h]
\centering{
\includegraphics[totalheight=1.5cm]{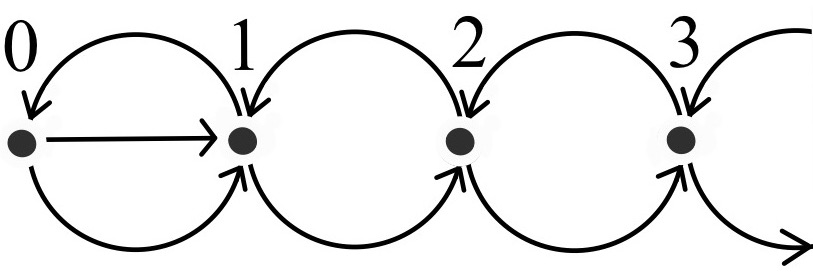}
\caption{The Cayley graph $\Gamma_{X,\{1\}}.$ }\label{1}
}
\end{figure}

\begin{example}\label{ex N and 0}
Let us take $X=\mathbb{N}\cup\{0\}$ with a 2-valued group operation $x*y=\bigl[x+y,|x-y|\bigr]$ (see \cite{B}). This is a 2-valued group with $0$ as the unit element and $inv(x)=x.$ It is clear that for any $x\in \mathbb{N}$ we have 
$$
x\in 1^{*x} = \underbrace{1 * 1 * \ldots * 1}_x,
$$ 
so $X$ is a 2-valued cyclic group generated by $1.$ Let us construct the Cayley graph of this group with respect to this generator (see  Figure~\ref{1}). 

We have $x*1=[x-1,x+1]$ for any $x>0,$ and $0*1=[1,1].$ Thus, vertices $x$ and $y$ are connected by an edge in $\Gamma_{X,\{1\}}$ iff $|x-y|=1.$ We now see that for this group,
$$\bigl|B(x,r)\bigr|=1+r+\min(x,r),$$
which  depends on $x.$ 
\end{example}

Still, we can show that the asymptotic properties of the growth function $\gamma_{\Gamma, x}$ do not depend on the choice of the element $x$ and the generating set $S.$

\begin{proposition} \label{asimpt}
Let $X$ be a $n$-valued group with finite generating sets $S$ and $S'.$ Let $y,y'\in X.$ Denote by $B(x,r)$ and by $B'(x,r)$ balls of radius $r$ centered at $x$ with respect to generating sets $S$ and $S'$ respectively. Then there is a constant $C > 0$ such that
$$|B(y, C^{-1} r)|\leqslant |B'(y',r)| \leqslant |B(y,C r)|.$$
\end{proposition}
\begin{proof}
Firstly, let us show that $B(e,r)$ grows at the same rate as $B(y,r)$. Since $S$ is a generating set of $X,$ both $y$ and its inverse $\bar{y}$ have well-defined lengths in $X$ with respect to $S.$ Since $y\in B\bigl(e,l(y)\bigr)$ and $e\in B\bigl(y,l(\bar{y})\bigr),$ we have
\begin{gather*}
|B(y,r)|\leqslant |B(e,r+l(y))|, \\
|B(e,r)|\leqslant |B(y,r+l(\bar{y}))|,
\end{gather*}
and it follows that for $l=1+\max\bigl(l(y), l(\bar{y})\bigr),$ and $r\geqslant 1$
$$|B(y,l^{-1} r)|\leqslant |B(e,r)| \leqslant  |B(y,l r)|.$$

Similarly, denote by $l(S')$ the maximum length of an element of $S'$ with respect to~$S,$ and by $l'(S)$ the maximum length of an element of $S$ with respect to~$S'.$ Denote $l=1+\max\bigl(l(S'),l'(S)\bigr)$ Now we have $S'\subset B\bigl(e,l(S')\bigr)$ and $S\subset B'\bigl(e,l'(S)\bigr),$ and
$$|B(e,l^{-1} r)|\leqslant |B'(e,r)| \leqslant |B(e, l r)|.$$

By the arbitrary choice of $y$ and $S,$ this concludes the proof.
\end{proof}

A famous result of 
 Gromov  \cite{G} says that any virtually nilpotent  group has a polynomial growth. It is natural to ask: is it true that any multivalued group which comes from a virtually nilpotent  group   has a polynomial growth? For coset groups an answer gives

\begin{theorem} \label{Gr}
Let $G$ be a finitely generated virtually nilpotent group, let $A$ be a finite subgroup of $\Aut(G)$ of cardinality $n$. Then the  $n$-valued coset group $X = (G, A)$   has a polynomial growth.
\end{theorem}

\begin{proof}
Let $G$ be generated by a  finite set  $S =S^{-1} =\{ s_1, s_2, \ldots, s_m \}$ and   $GA$ be  the semi-direct product of $G$ and $A$ that is a group which consists  of the  pairs $(g, a) \in G \times A$ with a multiplication 
$$
(g, a) (h, b) = (g h^{a^{-1}}, a b), ~~g, h \in G, a, b \in A.
$$
Here we use a notation $h^{a^{-1}} = a^{-1}(h)$ which means the image of $h$ under the action of the automorphism $a^{-1} \in \Aut(G)$.
Since $GA$ is a finite extension of $G$, it is virtually nilpotent. We consider the Cayley graph of $GA$ with respect to the set of generators $\{ (s_i, a) ~|~ s_i \in S,  a \in A, \}$. 
By Gromov's theorem the group $GA$ has a polynomial growth, i. e.  there is a natural numbers $k$ and $C$  such that for any  $r \geq 1$ holds  $|B(e, r)| \leq C r^k$.

Using Proposition \ref{asimpt} we can consider the balls $B_X(e, r)$ in the Cayley graph $\Gamma = \Gamma_{X,  \overline{S}}$, where $ \overline{S} = \{ \overline{s_1}, \overline{s_2}, \ldots, \overline{s_m} \}$ and  $\pi(s_i) = \overline{s_i}$. Since some $\overline{s_i}$ can be equal, we get $|\overline{S}| \leq |S| $.

Let us prove that $|B_X(e, r)| \leq |B_{GA}(e, r)|$ for $r = 1, 2, \ldots$. Induction by $r$. If $r=1$, then $B_{GA}(e, 1) = \{ e, (a, s_1), (a, s_2), \ldots, (a, s_m)~|~ a \in A \}$ 
and $|B_{GA}(e, 1)| \leq 1 + nm $.
 On the other hand,  
$$
B_X(e, 1) = \{ e, e * \overline{s_1}, e * \overline{s_2}, \ldots, e * \overline{s_m} \} = \{ e, \overline{s_1}, \overline{s_2}, \ldots, \overline{s_m} \}
$$
 and $|B_X(e, 1)| \leq 1 + m$. Hence, $|B_X(e, 1)| \leq |B_G(e, 1)|$.

Let us find $B_{GA}(e, 2)$. It is equal to $B_{GA}(e, 2) = B_{GA}(e, 1) \cup \{ (s_i, a) (s_j, b) ~|~ 1 \leq i, j \leq m,~~a, b  \in A \}$. Since, $G$ has polynomial growth, $|B_{GA}(e, 2)| \leq C 2^k$ for some $C$ and $k$.

On the other hand, 
$$
B_X(e, 2) = B_X(e, 1) \cup \{ \overline{s_i} * \overline{s_j} ~|~ 1 \leq i, j \leq m\}.
$$
By the definition of $n$-multiplication in $X$,
$$
\overline{s_i} * \overline{s_j} = \pi (\pi^{-1}(\overline{s_i} ) \cdot \pi^{-1}(\overline{s_j} )) = \pi (\{ s_i A \cdot s_j A~|~1 \leq i, j \leq m  \}).
$$
Since, $s_i A \cdot s_i A \subseteq \{ (s_i, a) (s_j, b) ~|~ 1 \leq i, j \leq m,~~a, b  \in A \}$ and $\pi$ is a projection, we get $|B_X(e, 2)| \leq |B_X(e, 2)|$.

Suppose that 
$$
|B_X(e, r)| \leq |B_{GA}(e, r)| \leq C r^k
$$
and consider the $r+1$-balls $B_X(e, r+1)$ and $B_{GA}(e, r+1)$. We have 
$$
B_X(e, r+1) = B_X(e, r) \cup \{ \alpha * \overline{s_j} ~|~ \alpha \in  B_X(e, r),  1 \leq  j \leq m\}.
$$
By the definition of $n$-multiplication in $X$,
$$
\alpha * \overline{s_j} = \pi (\pi^{-1}(\alpha) \cdot \pi^{-1}(\overline{s_j} )) = \pi (\{ \alpha A \cdot s_j A~|~1 \leq  j \leq m  \}).
$$
On the other hand,
$$
B_{GA}(e, r+1) = B_{GA}(e, r)  \cup \{ \beta (s_i, a)~|~ \beta \in B_{GA}(e, r), s_i \in S, a \in A \}
$$
and we have
$$
|B_{GA}(e, r+1)| \leq |B_{GA}(e, r)| + |B_{GA}(e, r)| nm \leq C (r+1)^k.
$$

Since, $\alpha A \cdot  s_j A \subseteq \{ \beta  (s_j, b) ~|~ 1 \leq  j \leq m,~~ b  \in A \}$ and $\pi$ is a projection, we get 
$$
|B_{X}(e, r+1)| \leq |B_{X}(e, r)| + |B_{X}(e, r)| nm.
$$ 
Hence,
$$
|B_X(e, r+1)| \leq |B_{GA}(e, r+1)|  \leq C (r+1)^k.
$$
\end{proof}

\begin{question}
Let $G$ be a group of polynomial growth, $H$ be its subgroup of cardinality~$n$. Is it true that the double coset group of the pair $(G, H)$ has a polynomial growth?
\end{question}


\section{$n$-valued dynamics} \label{dyn}

An $n$-valued dynamic $T$ on a space $Y$  is a map $T \colon Y \to (Y)^n$. For $y \in Y$ a function 
$\xi_y \colon \mathbb{N} \to \mathbb{N} $, where {$\xi_y(r)=\bigl|\Set\bigl(T^{\,r}(y)\bigr)\bigr|$} is called a growth function of $T$ in the point $y$. For any $n$-valued group $X$ and element $z\in X$ we can define  a $n$-valued dynamic  $T_z$ on $X$ by the rule
$$
T_z(y) = y * z,~~y \in X.
$$
Here we assume that
$$
T^2_z(y) = T_z(T_z(y)) = (y * z) * z\;, \ldots,\; T^{k}_z(y) = T_z(T^{k-1}_z(y)) = (T^{k-1}_z(y)) * z.
$$
In \cite{Bu} was formulated 

\begin{question} \label{Dyn}
Characterize $n$-valued dynamics $T$ with polynomial growth functions $\xi_y(r)$ for any $y \in Y$.
\end{question}

It is easy to see that in a $n$-valued group $X$ for any $z\in X$ the growth function $\xi_y$ which is defined by the dynamic $T_z$ satisfies the inequality
$$
\xi_y(r) \leq \gamma_{X,y}(r) ~\mbox{for all}~r \in \mathbb{N}. 
$$ 
As a corollary of Theorem \ref{Gr} we have

\begin{proposition} \label{poldyn}
If $X = (G, A)$, $A \leq \Aut G$, $|A| = n$, is a $n$-valued coset group, then any $n$-valued  dynamic $T_z \colon X \to (X)^n$ which is defined by an element   $z \in X$ has the  polynomial growth function $\xi_y(r)$ in  any $y \in X$.
\end{proposition}

For arbitrary group $G$ growth of $n$-dynamic  on coset group gives 

\begin{theorem}\label{Th growth function bounds}
Let $G$ be a group, $A\leq \Aut(G), |A|=n, g\in G.$ Then in the coset group $X=(G,A)$ for the element $z=gA$ and for any $y\in X$ the growth function $\xi_y(r),$ which is defined by the dynamic $T_z,$ satisfies the inequality:
$$
\frac{1}{n}|S^+(e, r)|\leq\xi_y(r)\leq |B^+(e, r)|,
$$ 
where $S^+(e, r)$ and $B^+(e, r)$ are the sphere and the ball respectively  in the monoid $M=\bigl\langle\{ g a  \mid a\in A\}\bigr\rangle\subset G.$
\end{theorem}
\begin{proof}
Let $y=hA$ for some $h\in G.$ We will prove by induction on $r$ that the map $\varphi\colon G\rightarrow X$ sending $p\in G$ to $hpA\in X$ maps $S^+(e,r)$ to the set $\Set (y*z^{*r})$. For $r=0$ we have $\varphi(e)=hA=y.$ 

Assume that the induction hypothesis holds for $r.$ Then, for any $p\in S^+(e,r+1),$ we have $p=q(g)a_1$  for some $a_1\in A$ and $q\in S^+(e,r).$ By the induction hypothesis, $qA\in \Set(y*z^{*r}).$ Then for any $a\in A$ we have 
$$(q(g)a)A\in \Set( qA*z)\subset\Set(y*z^{*(r+1)}),$$
which proves the induction hypothesis for $r+1.$

Left multiplication by $h$ is an automorphism  of $G,$ and the map sending $p\in G$ to its orbit under the action of $A$ can send no more than $n$ elements of $G$ to the same element of $X.$ Thus, we have $\frac{1}{n}|S^+(1,r)|\leq\xi_y(r).$

We will now prove by induction that $\Set(y*z^{*r})\subset \varphi\bigl(B^+(e,r)\bigr).$ For $r=0$  we have $\{y\}=\{h\}A=\varphi(e).$
Assume that the induction hypothesis holds for $r.$ Fix $t\in \Set(y*z^{*(r+1)}).$ There are $p\in G$ and $a\in A$ such that $t=(p(g)a)A$   and $pA \subseteq \Set(y*z^{*r}).$ By the induction hypothesis, $pA \in  \varphi\bigl(B^+(e,r)\bigr),$ i.e. there is $q\in B^+(e,r)$ such that $(hq)A=pA.$ Then there is $a'\in A$ such that $(p)a'=hq.$ Now we have
$$t=(p(g)a)A=\bigl((p(g)a)a'a'^{-1}\bigr)A=\bigl((p)a'(g)aa'\bigr)a'^{-1}A=\bigl((p)a'(g)aa'\bigr)A=\bigl(hq(g)aa'\bigr)A,$$
and, since $(g)aa'$ is one of the canonical generators of the monoid $M,$ we have $q(g)aa'\in B^+(e,r+1)$ and thus $t=\varphi(q(g)aa')\in \varphi\bigl(B^+(e,r+1)\bigr),$ which proves the induction hypothesis for $r+1.$

Since $\Set(y*z^{*r})\subset \varphi\bigl(B^+(e,r)\bigr),$ we have $\xi_y(r)\leq |B^+(e,r)|,$ which completes the proof.
\end{proof}

\bigskip

The cyclically presented group $G_n(w)$ is the group defined by the cyclic presentation
$$
\langle g_1, g_2, \ldots, g_n ~|~w, w a, \ldots, w a^{n-1} \rangle
$$
where $w = w(g_1, g_2, \ldots, g_n)$ is a word in the free group $F_n$ with generators $g_1, g_2, \ldots, g_n$ and $a \colon F_n \to F_n$ is the shift automorphism of $F_n$ given by $g_1a=g_{2},$ $g_2a=g_{3},$ $\ldots,$ $g_{n-1}a=g_{n}$ and $g_na=g_1$.

In \cite{BV}  Buchstaber  and Vesnin ask a question on growth functions of multivalued dynamics defined by cyclically presented groups. Let us show that Theorem \ref{Th growth function bounds} gives an answer to this question.
Let $G = G_n(w) =\langle g_1, \ldots, g_n\rangle$ be a cyclically presented group and suppose that the shift $a  \colon G \to G$   defines an automorphism of $G$.
Put $A=\langle a\rangle \leq \Aut(G)$. Then, the coset group $(G, A)$ is generated by a single element, that is the orbit of the generating set of $G$ under the action of $A$. By Theorem \ref{Th growth function bounds}, we have.

\begin{corollary} \label{cor4.4}
Let $G=\langle g_1, \ldots, g_n\rangle$ and $A=\langle a\rangle \leq \Aut(G)$ are defined above. Then  
 the $n$-valued dynamic  $T_z$, $z = g A$, $z\not=e$ of the coset group $(G, A)$ has an  exponential growth if and only if the group $G$ has an exponential growth. This dynamic $T_z$ has a polynomial growth if and only if the group $G$ has a polynomial growth.
\end{corollary}

There are cyclically presented groups with polynomial growth and with exponential growth, hence a dynamic of an element of a coset group $(G, A)$, where $G$ is a cyclically presented group, can have exponential or polynomial growth.

Theorem \ref{Th growth function bounds} allows us to construct a group $G$ with exponential growth and a coset group over $G$ without dynamics of exponential growth. 

Take in attention Proposition \ref{poldyn}, it is natural to formulate

\begin{conjecture} \label{con4.4}
Let $G$ be a group of exponential growth, $A \leq \Aut(G)$, $|A| = n$. Then any $n$-valued dynamic which is defined on the coset group $X = (G, A)$ by some element $z \in X$ has an exponential growth.
\end{conjecture}

The next example shows that this conjecture does not hold.

\begin{example} \label{ex}
Let $H$ be a finite group, $F$ be a free group of rank 2.  Put $G=H \times F.$ Now, let $A$ be a group of automorphisms of $G$ that acts trivially on $F.$ For any $g\in G,$ the subgroup $\bigl\langle\{g a\mid a\in A\}\bigr\rangle\subset G$ is isomorphic to $H_0 \times F_0$, where $H_0$ is finite and $F_0$ is either trivial or cyclic. No submonoid of $H_0 \times F_0$ has exponential growth. Therefore, no element of the coset group $X=(G,A)$ defines a dynamic with exponential growth.

In particular, let $H = \langle h \rangle$ be  a cyclic group of order 3 and $A = \langle a \rangle \leq \Aut (H \times F)$, where $a$ acts by the rule $(h, f) a = (h^{-1}, f)$ for any $f \in F$. It is easy to see that $A$ is the cyclic group of order 2  and the 2-valued dynamic $T_z$, which is defined by element $z = (h, e)$ does not have an exponential growth.
\end{example}

For an $n$-valued group $X$ and its element $x$ define $B^*(x,0)=S^*(x,0)=\emptyset,$ and for any $r\in \mathbb{N}$ define 
$$
B^*(x,r)=\bigcup\limits_{i=1}^r \Set(x^{*r}),~~S^*(x,r)=B^*(x,r)\setminus B^*(x,r-1).
$$
 It is clear that $S^*(x,r)\subset \Set(x^{*r})\subset B^*(x,r),$ and that $B^*(x,r)=\bigcup\limits_{i=1}^r S^*(x,i).$

\begin{lemma}\label{lm addition of spheres}
Let $X$ be an $n$-valued group, $x\in X.$

\textnormal{a)} If $S^*(x,r)=\emptyset$ for some $r\in\mathbb{N},$ then $S^*(x,r')=\emptyset$ for any $r'>r.$

\textnormal{b}) For $r_1,\ldots,r_k\in\mathbb{N}$ such that $S^*(x,r_i)\neq \emptyset,$ the following expression holds:
$$S^*(x,r_1+\ldots+r_k)\subset\Set\Bigl(S^*(x,r_1)*\ldots * S^*(x,r_k)\Bigr).$$
\end{lemma}
\begin{proof}
a) It suffices to show that $S^*(x,r+1)=\emptyset.$ If $y\in B^*(x,r+1),$ then $y=z*x$ for some $z\in\Set(x^{*r})\subset B^*(x,r).$ Since $S^*(x, r)=\emptyset$ and thus $B^*(x,r-1)=B^*(x,r),$ the element $z$ is a member of $\Set(x^{*r_0})\subset B^*(x,r_0)$ for some $r_0<r.$ Therefore, $y\in\Set(x^{*(r_0+1)})\subset B^*(x,r).$ We have shown that $$S^*(x,r+1)=B^*(x,r+1)\setminus B^*(x,r)=\emptyset.$$

b) Fix any $y\in S^*(x,r_1+\ldots+r_k).$ We have
$$y\in \Set(x^{*(r_1+\ldots+r_k)})=\Set(x^{*r_1}*\ldots*x^{*r_k})=\Set\Bigl(\Set(x^{*r_1})*\ldots*\Set(x^{*r_k})\Bigr),$$
and there is a set $\{z_1,\ldots,z_k\},$ such that $z_i\in\Set(x^{*r_i})$ for every $i\in\{1,\ldots,k\},$ and that $y\in \Set(z_1*\ldots*z_k).$

Suppose that for some $p\in\{1,\ldots,k\}$ there is $r_p'<r_p$ such that
$z_p\in \Set(x^{*r_p'}).$ In this case, $y\in \Set(x^{*(r_1+\ldots+r_{p-1}+r_p'+r_{p+1}+\ldots+r_k)}),$ and since $$r_1+\ldots+r_{p-1}+r_p'+r_{p+1}+\ldots+r_k<r_1+\ldots+r_k,$$ we have $y\in B^*(x,r_1+\ldots+r_k-1),$ which contradicts the choice of $y.$ Therefore, $z_i\in S^*(x,r_i)$ for every $i\in \{1,\ldots,k\}.$
\end{proof}

\begin{theorem} \label{t4.8}
Let $X$ be a $2$-valued group such that $inv(a)=a$ for any $a\in X.$ For any $x\in X$ the growth function $\xi_x(r),$ which is defined by the dynamic $T_x,$ satisfies the inequality $\xi_x(r)\leq r(r+1).$
\end{theorem}
\begin{proof}
Since for any $a\in X$ there is a unique element $b\in X$ such that $a*a=[e,b],$ we can inductively define the sequence $\{x_0,x_1,x_2,\ldots\}$ by assuming $x=x_0$ and $x_{k-1}*x_{k-1}=[e,x_k]$ for every $k\in \mathbb{N}.$

We will now explore the cardinality of the set $S^*(x,2^k)$. For $k=0$ we have $S^*(x,2^0)=\{x\}.$ For $k=1$ we have $$S^*(x,2^1)=\Set(x*x)\setminus\{x\}=\Set\bigl([e,x_1]\bigr)\setminus\{x\}.$$

Therefore, $|S^*(x,2^1)|\leq 2$ and $B^*(x,2^1)=\Set\bigl([e,x,x_1]\bigr).$

For $k=2$ we can use Lemma \ref{lm addition of spheres} to see that $S^*(x,2^2)$ is either empty or 
$$S^*(x,2^2)\subset \Bigl(S^*(x,2^1)*S^*(x,2^1)\Bigr)\subset\Set\bigl([e,x_1]*[e,x_1]\bigr).$$

Since $e*e=[e,e],$ $e*x_1=[x_1,x_1]$ and $x_1*x_1=[e,x_2],$ and 
$$\Set\bigl([e,x_1]\bigr)\subset B^*(x,2^1)\subset B^*(x,2^2-1),$$
we can conclude that the set $S^*(x,2^2)$ is either empty or equal to the one-element set~$\{x_2\}.$

Now assume that for some $k\geq 2,$ the set $S^*(x,2^k)$ is either empty or equal to $\{x_k\}.$ By lemma \ref{lm addition of spheres}, if $S^*(x,2^k)=\emptyset,$ then $S^*(x,2^{k+1})=\emptyset,$ and if $S^*(x,2^k)=\{x_k\},$ then 
$$S^*(x,2^{k+1})\subset\Set(x_k*x_k)=\Set\bigl([e,x_{k+1}]\bigr),$$
and since $e\in S^*(x,2^2)\subset B^*(x,2^{k+1}-1),$ we can conclude that the set $S^*(x,2^{k+1})$ is either empty or equal to the one-element set $\{x_{k+1}\}.$ We have now proven by induction that $|S^*(x,2^k)|\leq 1$ for $k\geq 2.$

Now take an arbitrary integer $r > 2.$ Assume that the set $S^*(x,r)$ is nonempty. Expand $r$ into powers of two: $r=\sum\limits_{i=0}^R 2^ir_i,$ where $r_i\in \{0,1\}$ and $R$ is the integer part of $\log_2r+1.$ By Lemma \ref{lm addition of spheres} we have 
\begin{equation}\label{eq decomposition of a sphere}
S^*(x,r)\subset\Set\bigl(S^*(x,2^{i_1})*\ldots*S^*(x,2^{i_p})\bigr),
\end{equation}
where $\{i_1,\ldots i_p\}$ is the nonempty set of all indexes $0\leq i\leq R$ such that $r_i=1.$ Note that the $n$-valued product in the expression \ref{eq decomposition of a sphere} is the product of no more than $R+1$ multi-sets, with possibly one of those multi-sets being a two-element multi-set and all the others being one-element sets. Therefore, we have
$$|S^*(x,r)|\leq 2^{R+2}\leq 2^{\log_2 r + 3} < 2r.$$
Note that the above inequality still holds if $r\leq 2$ or if $S^*(x,r)=\emptyset.$

Now since $$\Set(x^{*r})\subset B^*(x,r)=\bigcup\limits_{i=1}^r S^*(x,i),$$
we have $\xi_x(r)\leq \sum\limits_{i=1}^r 2i = r(r+1).$
\end{proof}
\color{black}

\bigskip



\begin{ack}
This work has been supported by the grant the Russian Science Foundation, RSF 24-21-00102, https://rscf.ru/project/24-21-00102/.
\end{ack}

\end{document}